# Security Analysis and Fault Detection Against Stealthy Replay Attacks

Amirreza Zaman, Behrouz Safarinejadian

School of Electrical and Electronic Engineering, Shiraz University of Technology, Modarres Blvd., P.O. Box: 71555-313, Shiraz, Iran. Email: amirrezazaman@gmail.com

**Abstract**

This paper investigates the security issue of the data replay attacks on the control systems. The attacker is assumed to interfere with the control system's process in a steady-state case. The problem is presented as the standard way to attack, which is storing measurements and replaying them in further times to the system. The controller is assumed to be the LQG controller. The main novelty in this paper can be stated as proposing a different attack detection criterion by using the Kullback-Leibler divergence method to cover more general control system problems with these attacks and with higher-order dynamics. Also, there exists a packet-dropout feature in transmitting the data as another contribution of the paper. Formulations and numerical simulations prove the effectiveness of the newly proposed attack detection procedure by having a quick response to occurred attacks. Although, in previous approaches, the trade-off between attack detection delay or LQG performance was significant, in this approach it is proved that the difference in this trade-off is not considered in early moments when the attack happens since the attack detection rate is rapid and thus, these attacks can be stopped with defense strategies in the first moments with the proposed attack detection criterion.

Keywords: Cyber-physical system security, Kullback-Leibler divergence, LQG controller, Replay attacks,

**1. Introduction**

To totalize computations, control procedures, communication networks, and physical processes, cyber-physical systems (CPSs) are used as the next generation of engineered systems [1]. So, cyber-physical systems are essential nowadays since based on the data transfer procedure between different technologies and components, and CPS can be threatened by different malicious attackers to change the components' data or steal them which may lead to intense irrecoverable results on the global economy, social security or even human lives irreparable damages [2].

During the past decade, because of the rising usage of CPS, attack strategies and defense mechanism designings are increased significantly as well. Attackers in the system consideration can be determined as malicious agents in Wireless Sensor Networks. From the attacker's side, the goal is to degrade the system performance considerably under a condition that attackers remain undetected from the system fault detector at each time step $k$ [3]. One type of attacks is denial-of-service attacks. These kinds of attacks try to block the communication network to barricade the accurate access to the system components.

It should be noted that the jamming process needs a lot of transmission power in networks. Since the available energy for data transfer is limited, it can often be impossible to define and analyze such attacks. Besides, attack models based on jammer for attackers which considered to be resource-constrained were investigated [4,5]. Furthermore, some studies have been done about obtaining optimal data-forwarding programming against remote state estimation with game-theoretic formulation by [6,7].

Another type of cyber attacks can be organized as deception attacks. In these types of attacks, the attacker's goal is to threaten the authenticity of the sensors and actuators' data by injecting false data among them. One control plan against vague deception attacks of nonlinear time-delay CPSs is studied in [8], in which an adaptive, resilient, dynamic surface control using the neural-network scheme is formulated for these attacks on both sensor and actuator. The problem of designing a distributed filter for linear discrete-time networked control systems against deception attacks and bounded disturbances is investigated in [9] using a Round-Robin-type protocol to monitor communication between filters because of the communication resource constraint. In [10], an online attack strategy is designed and proved that the procedure can degrade the estimation quality by analyzing the proposed cost deviation, which showed the difference between the estimation quality with and without the proposed attack strategy. The attack scheduling of deception attacks for discrete-time systems with attack detection is studied in [11] particularly for a class of Kalman filters with $\chi^2$ detectors by assuming different types of attack scenarios (i.e., consecutive deception attacks or randomly launched deception attacks).

Furthermore, in [12], a distributed state estimator with an event-triggered scheme is formulated to defend against false data injection attacks in Wireless Sensor Networks. An approach [13] has been made recently to investigate security problems against deception attacks. In this approach [13], instead of defining a unique cyber-attack model, authors focused on analyzing the system's response during false data injection attacks. Finally, they proposed the necessary and sufficient conditions under which the attacker could perform each kind of attack without being detected. Some other classes of attacks can be defined as replay attacks. According to its name, a replay attack aims to inject fake and unreal input control signals. Also, they attempt to replay to the past sensory data to keep stealthiness altogether. Some studies have been done about obtaining suitable conditions and countermeasure of replay attacks by [14,15] for linear Gaussian systems. To analyse the interchange between control efficiency and system security, some studies were done under a random game frame by [16]. In the case of electric power grids, work was done to inject false data into the system versus state estimation at the remote estimation center by [17]. Analyzing malicious data injection attacks on remote estimation performance was investigated in [18]. Some remarkable works have been done about studying trade-off between attack sneakiness and system performance set back for control signal infusion attack for first-order systems under $\epsilon$- weakly marginal stealthiness measurement scale by [19,20], and for higher-order dynamic systems with linear attack strategies by [21]. In addition to that, developing integrity attacks and analyzing results for secure state estimation was studied by [22,23].

It is worth noticing that a worst-case attacker in control CPSs can be designed optimally, and therefore, it can degrade the control system's performance disruptively. This matter has gained much attention to consider later, and a few approaches are devoting to this matter to struct effective defense strategies. So, obtaining useful intercepted-data detectors in the existence of smart attackers is still unverified.

In this paper, the security issue of the data replay attack on control systems is investigated. The replay attacker is considered to interfere with the control system's operation in a steady state. The attacker tries to read the sensors' data and save their readings for a specific amount of time, and then, the attacker injects this stored data to the control system. Thus, the input data to the control system will not be the authentic real-time observed data, and the control system is not aware of this issue that can happen at any time. Therefore, by providing the wrong input control to the actuator based on the false data, the control system's performance will be degraded quickly. This kind of attack is the most common attack occurring in the control systems. Also, the implementation of this attack is straightforward since the replay attacker just needs to be aware of the system's steady-state behavior and not the whole system's dynamics. The control system in this approach is assumed to be a discrete-time linear time-invariant Gaussian control system using an infinite time horizon Linear Quadratic Gaussian (LQG) controller. As stated before, the outcomes on replay attacks' detectability have not been defined in an integrated security scheme, and there represents not to settle a consent on the attack sneakiness' measures. Previous work regarding detecting a replay attack [24] applied a specific detection classic $\chi^2$ detector scheme as a detection criterion for a particular group of replay attacks. However, [24] did not study the cases when attack policies are not rigorously stealthy to the attack's detector. Different from the previous approach, in our work, the control system is equipped with the Kullback-Leibler (K-L) divergence false-data detector as a measure of attack sneakiness.

Moreover, in this paper, the packet-dropout feature of the communication network is considered to derive formulas for more general cases as another case of novelty. The presented method of replay attack's detection in this work is proved to be independent of any predefined detection frameworks. Therefore, systems with general dynamics and higher-orders can be considered to provide security schemes for them based on our approach for replay attacks. Besides, it is proved and illustrated that the detection response with the K-L divergence criterion and especially in this work to replay attacks, becomes so much faster than any other presented statistical detection methods.

The rest of this approach is organized as the following: In Section 2, the model of the control process is presented, and then, the state estimation procedure by considering the packet-dropout feature is analyzed. Also, in Section 2, the problem is formulated by applying Kalman filter equations and optimal LQG controller with the K-L divergence fault detection criterion. In Section 3, the straightforward model of the malicious replay attack is investigated, and its influence on the control system is discussed. In Section 4, the fault detection scheduling scheme is presented in the existence of a replay attack, and its usefulness is demonstrated in Section 5. Eventually, in Section 6, the conclusion of this study is presented.

## 2. Process' analytical model and state estimation procedure

In this section, the proposed problem is formulated by applying Kalman filter equations, the conventional LQG controller, and the K-L divergence fault detector for replay attacks.

Assume the general linear time-invariant (LTI) control system is described as follows:

$$x(k + 1) = Ax(k) + Bu(k) + w(k)$$

$$y(k) = Cx(k) + v(k)$$

where $x(k) \in \mathbb{R}^n$ is defined as the general vector of the state variables at each time step $k$, $w(k) \in \mathbb{R}^n$ is the noise of the process at time step $k$. Also, $y(k) \in \mathbb{R}^m$ is measurements' vector with the measurement noise, which is defined as $v(k)$ at time step $k$. It is assumed that process and measurement noises are independent of each other and denoted with $w(k) \sim N(0, W)$ and $v(k) \sim N(0, V)$.

### 2.1. Kalman filter equations by considering packet-dropout feature in a communication network

At every step $k$, the sensor transmits its measured data to a remote estimator. To estimate the system state in the remote estimator in a case that there is a packet dropout feature in the communication network [25], a Kalman filter is used to estimate the received data using the following equations:

$$\hat{x}^-(0) = \bar{x}(0) \,, P^-(0) = \Psi$$

$$\hat{x}^-(k) = A\,\hat{x}(k - 1) + Bu(k)$$

$$P^-(k) = AP(k - 1)A^T + W$$

$$K(k) = P^-(k)C^T(CP^-(k)C^T + V)^{-1}$$

$$\hat{x}(k) = \hat{x}^-(k) + \beta(k)K(k)(y(k) - C\hat{x}^-(k))$$

$$P(k) = P^-(k) - \beta(k)K(k)CP^-(k)$$

It is considered that $\beta(k)$ which is the packet-dropout coefficient, is known to the attacker and the system at each time step $k$. In this approach, it is considered that $\beta(k)$ has a fixed value during the time. So, $P(k)$ and $K(k)$ will converge to steady-state values as follows :

$$P = \lim_{k \to \infty} P^-(k) \,, \text{ and } K = PC^T(CPC^T + V)^{-1}$$

Because of the assumption that any control systems run for a long time horizon in studies and by considering a fixed packet-dropout value, it is assumed that the control system in our approach can be run at a steady state from the starting point. In this case, the initial condition for the covariance matrix is assumed to be $\Psi = P$. With this assumption, the Kalman filter transforms to a fixed-value standard gain estimator, therefore:

$$x^-(0) = \bar{x}(0) \,, \hat{x}^-(k) = A\hat{x}(k - 1) + Bu(k)$$

$$\hat{x}(k) = \hat{x}^-(k) + \beta K(y(k) - C\hat{x}^-(k)),$$

Where $\beta$ is the fixed packet-dropout coefficient of the communication channel.

### 2.2. Optimal control formulation using Linear Quadratic Gaussian ( LQG) controller

Based on the obtained state estimation $\hat{x}(k)$ An objective function can be defined. Then, the goal of the LQG controller is to minimize the defined objective function as follows.

$$\lambda = \min \lim_{t \to \infty} E \frac{1}{t} [\sum_{k=0}^{t-1} (x^T(k)Fx(k) + u^T(k)Gu(k))],$$

where $F$ and $G$ are defined as positive semi-definite matrices and $u(k)$ is assumed to be measurable regards to $y(0), \ldots, y(k)$. So, $u(k)$ can be obtained as a function of the above measurements. The answer to the defined minimization problem can be derived as a steady-gain controller like the following formulation :

$$u(k) = u^{opt}(k) = -(B^T RB + G)^{-1} B^T RA\hat{x}(k),$$

where $u^{opt}(k)$ is the obtained optimal control input and $R$ can be attained using the following Riccati equation:

$$R = A^T RA + F - A^T RB(B^T RB + G)^{-1} B^T RA .$$

Let define $M \triangleq -(B^T RB + G)^{-1} B^T RA$, then $u^{opt}(k) = M\hat{x}(k)$.

Now, the determined objective function by considering optimal estimator and the LQG controller in this approach would be rewritten as :

$$\lambda = Trace(RW) + Trace[(A^T RA + F - R)(P - \beta KCP)]$$

### 2.3. Intercepted-data detector

To monitor and control the system behavior at the remote estimator, we need to implement an intercepted-data detector to find out whether there are cyber-attacks or not. There are some ways to detect attacks in a cyber-physical system, such as $\chi^2$ detector or $K-L$ divergence between two probability distributions or game-theoric methods [6,7,17,21]. In this paper, we are going to use $K-L$ divergence method [19] to detect attacks in the system. The Kullback – Leibler divergence ($K-L$ divergence) is a non-negative distance measurement criterion between two probability distributions, which is practical to use in cyberattacks' tracing theories [19, 26].

**Definition 1. Kullback – Leibler Divergence.** This criterion expresses that if there are two stochastic sequences $x(k)$ and $y(k)$ with corresponding joint probability density functions $f_{x(k)}$ and $f_{y(k)}$, respectively, then, the Kullback – Leibler divergence between $x(k)$ and $y(k)$ will be defined as given below:

$$D\big(x(k) \parallel y(k)\big) = \int_{\{\mu(k) \mid f_{x(k)}(\mu(k)) > 0\}} \log \frac{f_{x(k)}(\mu(k))}{f_{y(k)}(\mu(k))} f_{x(k)}(\mu(k)) d\mu(k)$$

According to the above $K-L$ divergence equation, if $f_{x(k)} = f_{y(k)}$ then $\big(x(k) \parallel y(k)\big) = 0$. It is proved that $K-L$ divergence criterion is not symmetric in general, thus $\big(x(k) \parallel y(k)\big) \neq D\big(y(k) \parallel x(k)\big)$.

### 3. Linear replay attack strategy derivation against operating control system

In this section, to derive formulations, it is considered that an adverse third party tries to intrude on the presented control system. The replay attack model in this approach is defined as computer security problems. Also, the feasibility of these kinds of attacks on the control system is investigated. It is proved that the analysis of this work can be generalized to other types of control systems with higher orders. The attacker can inject the control input $u_a(k)$ at any time. It is worth noticing that the identification process of the underlying dynamic model of the control system for attackers can be hard in general, and not all the attackers are such powerful to detect systems' models. Therefore, this paper focuses on a straightforward attack strategy, which is much easy to implement. Besides, since the energy-consuming limitation exists in reality for both system and attacker, the goal of designing attack or defense strategies in CPSs is to develop with the lowest level of energy consumption for both parties.

In the existence of attack in control systems, to implement any counter-attack strategies, at first, the attacker should be detected as fast as it can. Thus, in the attack stage, the intercepted-data detector should work correctly. As mentioned, the attacker aims to inject fake stored measurements from any specific time step $k$ to the control system. To have an integrated framework for analyzing these kinds of attacks, the stored measures by the attacker can be considered as the output of a virtual network as given in the following:

$$\tilde{x}(k+1) = A\tilde{x}(k) + B\tilde{u}(k) + w(k)$$

$$\tilde{y}(k) = C\tilde{x}(k) + v(k)$$

It is assumed that $w(k)$ and $v(k)$ are independent of the attack. Therefore, they are the same as the ones described for the real system.

Also, the virtual system can be introduced as follows.

$\hat{\tilde{x}}^-(k) = A\hat{\tilde{x}}(k-1) + B\tilde{u}(k),$

$\hat{\tilde{x}}(k) = \hat{\tilde{x}}^-(k) + \beta K\left(\tilde{y}(k) - C\hat{\tilde{x}}^-(k)\right)$, and

$\tilde{u}(k) = M\hat{\tilde{x}}^-(k)$

with the specified initial conditions $\tilde{x}(0)$ and $\tilde{x}^-(0)$. If it is considered that the attacker can learn the system's behavior at each time step $k$, then it is evident that the attacker would run the virtual network to maximize the attack's influence on the system's performance. To investigate a replay attack, consider that the attacker can store the previous system's measurements. So, the mentioned virtual system would be just the time-shifted edition of the actual order. Let define T as the time shift. Then, it is easy to conclude that $\tilde{x}(k) = x(T+k)$ and $\tilde{x}^-(k) = x^-(T+k)$.

Assume that the system is under a replay attack and the system's defender is applying the K-L divergence criteria between the actual system's innovation sequence $z(k) = y(k) - C\hat{x}(k)$ and the input innovation sequence $\tilde{z}(k) = \tilde{y}(k) - C\hat{\tilde{x}}^-(k)$ which is sent via the virtual system to detect attack interference in the system. With this consideration, the Kalman filter equations can be rewritten as given below in a recursive way :

$\hat{x}(k) = A\hat{x}(k-1) + Bu(k) = (A + BM)\hat{x}(k-1) = (A + BM)\left[\hat{x}^-(k) + \beta K(\tilde{y}(k) - C\hat{x}^-(k))\right]$
$= (A + BM)(I - \beta KC)\hat{x}^-(k) + (A + BM)\beta K\tilde{y}(k)$

Furthermore, for the virtual system, a similar equation for $\hat{\tilde{x}}^-(k)$ can be derived as follows:

$\hat{\tilde{x}}(k) = (A + BM)(I - \beta KC)\hat{\tilde{x}}^-(k) + (A + BM)\beta K\tilde{y}(k) \ .$

For simplicity, it is assumed that the time of the attack starts as time 0. Define $\Lambda$ as $\Lambda \triangleq (A + BM)(I - \beta KC)$, then

$\hat{x}^-(k) - \hat{\tilde{x}}^-(k) = \Lambda^k \left(\hat{x}^-(0) - \hat{\tilde{x}}^-(0)\right).$

Let define $\hat{x}^-(0) - \hat{\tilde{x}}^-(0) \triangleq \eta$ , now the residue from the previous equations can be rewritten as follows.

$\tilde{y}(k) - C\hat{x}^-(k) = \left(\tilde{y}(k) - C\hat{\tilde{x}}^-(k)\right) - C\Lambda^k\eta$

According to the virtual system's explanation, it can be inferred that the residue $\tilde{y}(k) - C\hat{\tilde{x}}^-(k)$ has a similar distribution to the residue $y(k) - C\hat{x}^-(k)$. Hereupon, if $\Lambda$ is considered to be stable during the time horizon, the second term of the mentioned equation converges to zero. Therefore, the residue $\tilde{y}(k) - C\hat{x}^-(k)$ has the same distribution as $z(k) = y(k) - C\hat{x}^-(k)$. So, the attack innovation sequence $\tilde{z}(k) = \tilde{y}(k) - C\hat{\tilde{x}}^-(k)$ can be rewritten as $\tilde{z}(k) = \tilde{y}(k) - C\hat{x}^-(k)$. Moreover, because of the distribution similarity between innovation sequences $z(k)$ and $\tilde{z}(k)$, the calculated detection rate given by the K-L divergence intercepted-data detector will consist of the same false alarm detection rate as before. So, the proposed intercepted-data detector is ineffectual.

Now, the second case of $\Lambda$ is investigated. If $\Lambda$ is assumed to be unstable, the attacker is not able to carry on replaying the measurement $\tilde{y}(k)$ over the time horizon of which the intercepted-data detection criteria will transform into an unbounded rule quickly. In this case, it can be concluded that the system is resilient to this kind of attack, as the defined intruder detection method can detect the attack's existence. Besides, with derived formulas and calculations, it is proved that the obtained outcomes from the implementation on a specific group of the estimator, controller, and intercepted-data detector for our approach, can be generalized to any virtual systems. Moreover, the outcomes from the studied K-L divergence intrusion detection method in this paper can be applied to any nonlinear systems and higher-order systems since other detection methods like $\chi^2$ detector can only be used for first-order and linear systems, but with using the K-L divergence method, this issue can be solved.

**Theorem 1.** The K-L divergence criterion between two innovation sequences $z(k)$ and $\tilde{z}(k)$ can e obtained as follows:

$$D\big(\tilde{z}(k) \parallel z(k)\big) \cong \frac{1}{2} Trace(cov(z(k))^{-1} cov(\tilde{z}(k))) - \frac{m}{2} + \frac{1}{2}\log\frac{|cov(z(k))|}{|cov(\tilde{z}(k))|}$$

Where $m$ is the dimension of $z(k)$ as general.

*Proof.* The proof of this Theorem is stated in [21].

## 4. Fault detection scheduling of deception replay attack

As mentioned in previous sections, it is necessary to consider designing detection and defense strategies against cyberattacks in vulnerable control systems, especially replay attacks.

To derive the attack detection strategy, it is assumed that $\Lambda$ is stable. According to the control systems with the LQG controller, the weak point of the LQG controller and Kalman filter is that they operate based on a fixed control gain, or a control gain that converges soon. Thus, we can conclude that the control approach, in this case, would be static from some viewpoints. As previous methods regarding the detection of replay attacks, the controller is assumed to be as the following form:

$$u(k) = u^{opt}(k) + \Delta u(k)$$

where $u^{opt}(k)$ is the output of the LQG controller, and $\Delta u(k)$s are collected from an i.i.d Gaussian distribution with zero mean and covariance value of $\tau$. It should be noticed that selected $\Delta u(k)$s are independent of the optimal control input $u^{opt}(k)$. The whole system's diagram can be presented in Figure 1.

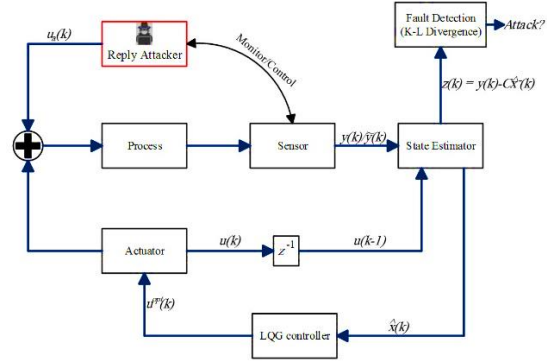

Figure 1. Control system's diagram

The purpose of adding the signal $\Delta u(k)$ is to verify the correct performance of the system. As mentioned before, $\Delta u(k)$ at each time step $k$ should be a zero-mean signal to prevent any bias to $x(k)$. It can be inferred that without any attacks, the LQG controller will not be optimal anymore. So, to detect any attacks, the controller's performance should be optimized, too. In the given Theorem [24], the loss of the used LQG controller is distinguished when the signal $\Delta u(k)$ exists.

**Theorem 2.** [24] By considering the existence of the signal $\Delta u(k)$ in the system, the LQG controller's performance can be stated as follows.

$$\tilde{\lambda} = \lambda + Trace[\,(G + B^T RB)\tau\,]$$

The proof of this Theorem is given in [24].

Now, different from the previous work, we consider the K-L divergence intercepted-data detector to consist of systems with more general features and higher-order dynamics after augmenting the accidental control signal. The following stated Theorem proves the usefulness of the K-L divergence intercepted-data detector under the altered control design.

**Theorem 3.** In the absence of any replay attacks or the case of occurring replay attacks without any new augmented detection strategies,

$$D\big(\tilde{z}(k) \| z(k)\big) = 0.$$

And under a replay attack,

$$D\big(\tilde{z}(k) \| z(k)\big) \cong 2 - m + Trace(\Sigma^{-1} 2C\Omega C^T) + \log|\Sigma| - \log|\Sigma + 2C\Omega C^T|,$$

where $m$ is the dimension of $z(k)$ as general, and $\Omega$ can be obtained by solving the given Lyapunov equation

$$\Omega - B\tau B^T = \Lambda\Omega\Lambda^T$$

And $\Sigma$ is the covariance of the innovation sequence of the real system and equals to $\Sigma = CPC^T + V$.

*Proof.* For the first equation, based on the fact that when there is no replay attack in the system or the case that there is a replay attack without any new augmented detection strategies, $\tilde{z}(k) = z(T + k)$, then $cov(\tilde{z}(k)) = cov(z(k))$. So, it can be inferred that:

$$D\big(\tilde{z}(k) \| z(k)\big) = 0.$$

This conclusion is an evident fact as a result of Gibb's inequality for the K-L divergence criterion.

For proving the second part of the Theorem 3, rewrite $\hat{x}(k)$ based on the augmented control input $\Delta u(k)$, as presented as a detection strategy as follows.

$$\hat{x}(k) = \Lambda \hat{x}^-(k) + (A + BM)\beta K \tilde{y}(k) + B\Delta u(k).$$

Therefore,

$$\hat{x}^-(k) - \hat{\tilde{x}}^-(k) = \Lambda^k \eta + \sum_{j=0}^{k-1} \Lambda^{k-j-1} B\big(\Delta u(j) - \Delta \tilde{u}(j)\big),$$

where $\eta = \hat{x}^-(0) - \hat{\tilde{x}}^-(0)$.

Now, $\tilde{y}(k) - C\hat{x}^-(k)$ can be rewritten as :

$$\tilde{y}(k) - C\hat{x}^-(k) = \tilde{y}(k) - C\hat{\tilde{x}}^-(k) - C\Lambda^k \eta - C\sum_{j=0}^{k-1} \Lambda^{k-j-1} B\big(\Delta u(j) - \Delta \tilde{u}(j)\big)$$

The first term of the above equation has just the same distribution as $y(k) - C\hat{x}^-(k)$. Also, by considering the stability of $\Lambda$, $C\Lambda^k \eta$ will be stable and converges to zero. Moreover, it is assumed that $\Delta u(j)$ is independent of the virtual system's dynamic and for the mentioned virtual system, $\tilde{y}(k) - C\hat{\tilde{x}}^-(k)$ will be independent of $\Delta \tilde{u}(j)$. Thus, the steady-state value of the virtual system's innovation sequence's covariance can be obtained as given in the following:

$$\varphi = \lim_{k\to\infty} cov\big(\tilde{y}(k) - C\hat{x}^-(k)\big) = \lim_{k\to\infty} cov\big(\tilde{y}(k) - C\hat{\tilde{x}}^-(k)\big) + \sum_{j=0}^{\infty} cov(C\Lambda^j B\Delta u(j)) + \sum_{j=0}^{\infty} cov(C\Lambda^j B\Delta \tilde{u}(j)) = \Sigma + 2\sum_{j=0}^{\infty} C\Lambda^j B\tau B^T (\Lambda^j)^T C^T,$$

where $\Sigma = CPC^T + V$. Let define $\Omega$ as follows.

$\Omega = \sum_{j=0}^{\infty} \Lambda^j B\tau B^T (\Lambda^j)^T$. By this definition, $\varphi$ is obtained as follows.

$$\varphi = \lim_{k\to\infty} cov\big(\tilde{y}(k) - C\hat{x}^-(k)\big) = \Sigma + 2C\Omega C^T.$$

Now, based on Theorem 1, the K-L divergence criterion between innovation sequences $\tilde{z}(k)$ and $z(k)$ can be obtained as follows.

$$D\big(\tilde{z}(k) \parallel z(k)\big) \cong \frac{1}{2} Trace(\Sigma^{-1}\varphi) - \frac{m}{2} + \frac{1}{2}\log\frac{|\Sigma|}{|\varphi|} \cong \frac{1}{2} Trace(\Sigma^{-1}[\Sigma + 2C\Omega C^T]) - \frac{m}{2} + \frac{1}{2}\log\frac{|\Sigma|}{|\Sigma + 2C\Omega C^T|}$$
$$\cong Trace(\Sigma^{-1}\Sigma + \Sigma^{-1}2C\Omega C^T) - m + \log|\Sigma| - \log|\Sigma + 2C\Omega C^T| \cong 2 - m + Trace(\Sigma^{-1}2C\Omega C^T) + \log|\Sigma| - \log|\Sigma + 2C\Omega C^T|$$

where $m$ is the dimension of $z(k)$ as general and $\Omega = \sum_{j=0}^{\infty} \Lambda^j B \tau B^T (\Lambda^j)^T$ and

$$\Sigma = CPC^T + V. \qquad \blacksquare$$

Moreover, one feasible way to design $\tau$ appropriately is to minimize equation $\tilde{\lambda} - \lambda = Trace[\,(G + B^T RB)\tau\,]$, while trying to maximize the value of $D\big(\tilde{z}(k) \parallel z(k)\big)$, which means maximizing the term $Trace(\Sigma^{-1}2C\Omega C^T)$.

## 5. Numerical Results

In this section, numerical simulations are provided to illustrate the analytical outcomes of this approach. Consider a comprehensive control system with two states described in section 2 with the following information :

$M$ is chosen as $M = -0.6180$. One candidate for selecting $\tau$ can be the steady-state covariance value of the system's state estimation without any attacks so that it can be chosen as $\tau = \begin{bmatrix} 1.1721 & 0.3146 \\ 0.3146 & 1.0229 \end{bmatrix}$. Although it is proved that the standard Kalman filter converges to its steady-state values very quickly, but the simulations in this paper are done for several iterations to monitor the system's behavior with or without a replay attack precisely. Also, it is considered that $F = G = I$. Without any control precautions against replay attacks, the system is defenseless to these attacks at first. It is also assumed that the attacker tries to begin replaying to the first ten state estimations at the filter's 10th iteration and keeps on attacking.

As it is proved in previous sections, especially in Theorem 3, the mentioned K-L divergence intercepted-data detector would have zero value under two cases. The first case, when there is no attack and second case when there is a replay attack without any defense strategies. The calculated value of the K-L divergence method under these two circumstances is depicted in Figure 2. Thus, it is observable that this detector is useless without taking into account any new defense strategies against a replay attack.

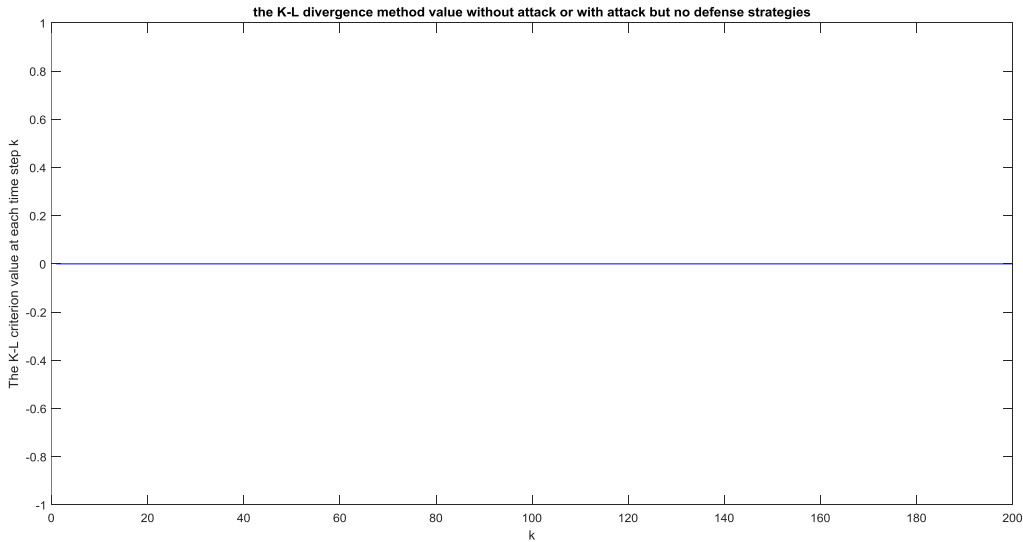

Figure 2

In Figure 3, the trace cost of the state estimation covariance of the real system is demonstrated as a replay attack's degrading measure when there is no attack and also when there is a replay attack to diverge estimation's process. As it is illustrated in Figure 3, the trace value of the system's covariance matrix under no attacks converges to the cost of 2.195. Still, when an attack occurs, the state estimation process is affected by the attacker, and the trace value is not going to converge during the time horizon. So, without a new detection strategy, the attacker can degrade the control system's performance while being stealthy in the system.

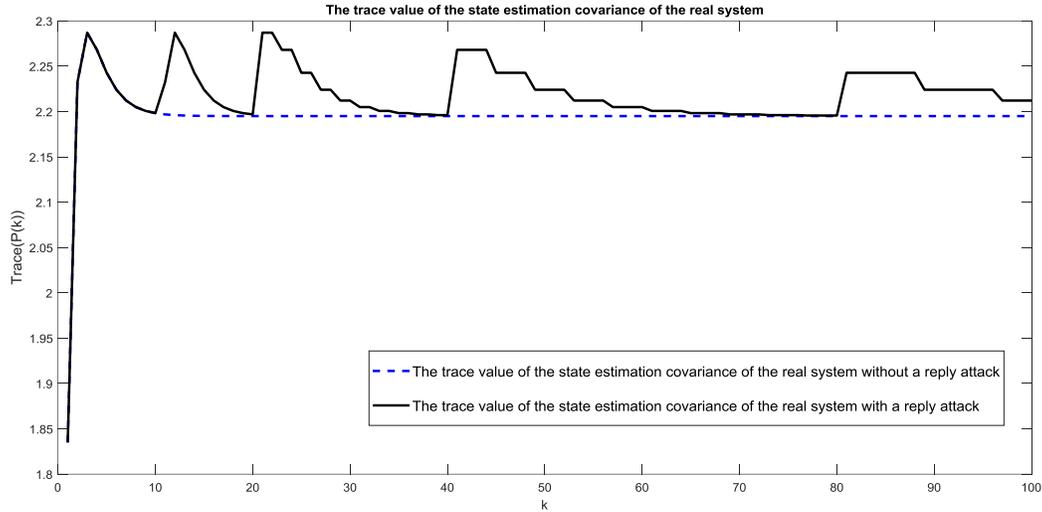

Figure 3

To compensate for the replay attack's influence on the system's performance, the newly presented control strategy is implemented, and the original value of the K-L divergence intercepted-data detector is calculated again based on Theorem 3. Figure 4 illustrates the computed value of the detection method based on the proposed control strategy. As it is shown in Figure 4, the K-L divergence attack detection method, in this case, has a specific extremum value in steady-state, which in this case, is 2.

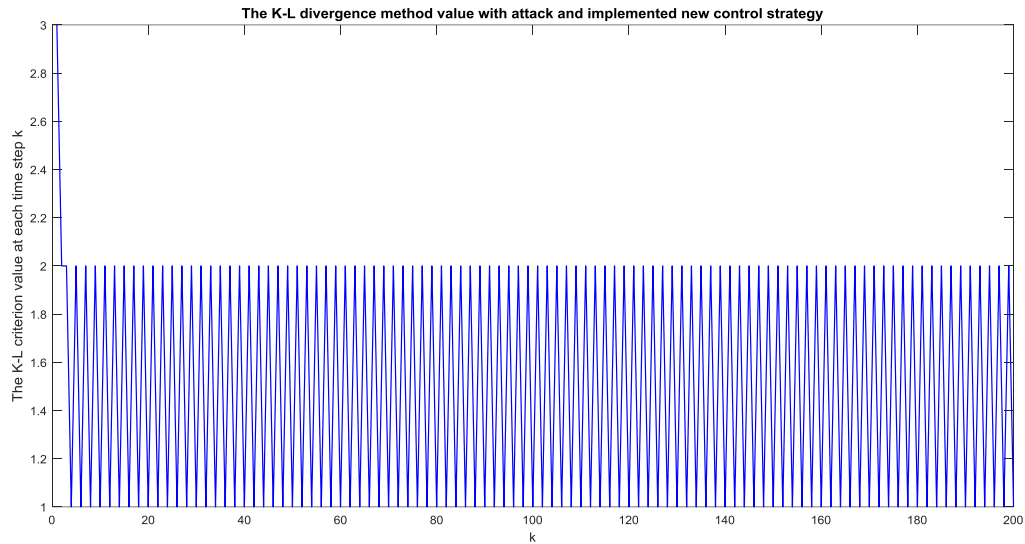

Figure 4

Thus, by considering a new control strategy, the system's estimation performance can be controlled, and also the replay attack can be detected very quickly. Besides, it is worth noticing that based on the previous approaches which used $\chi^2$ detector [ g], the response of the $\chi^2$ indicator depends on its defined window size. But it is proved and shown in this approach that since the K-L divergence criterion is independent of the system's dynamic or any detection factors. Also, the detection response with the K-L divergence criterion and especially in this work to replay attacks, becomes so much faster than any other statistical detection methods. Thus, the operator does not need to choose any detection window sizes to detect replay attacks suitably. Besides, the loss of the LQG performance after injecting signal $\Delta u(k)$ is illustrated in Figure 5. We can also conclude from Figure 5 that the control performance of the real system should be changed in a steady-state, and then it will not be optimal as the no-attack case. Also, since the attack happens in $10^{th}$ iteration, the signal $\Delta u(k)$ is injected from this iteration to detect an attack, and this matter can be observed in Figure 5.

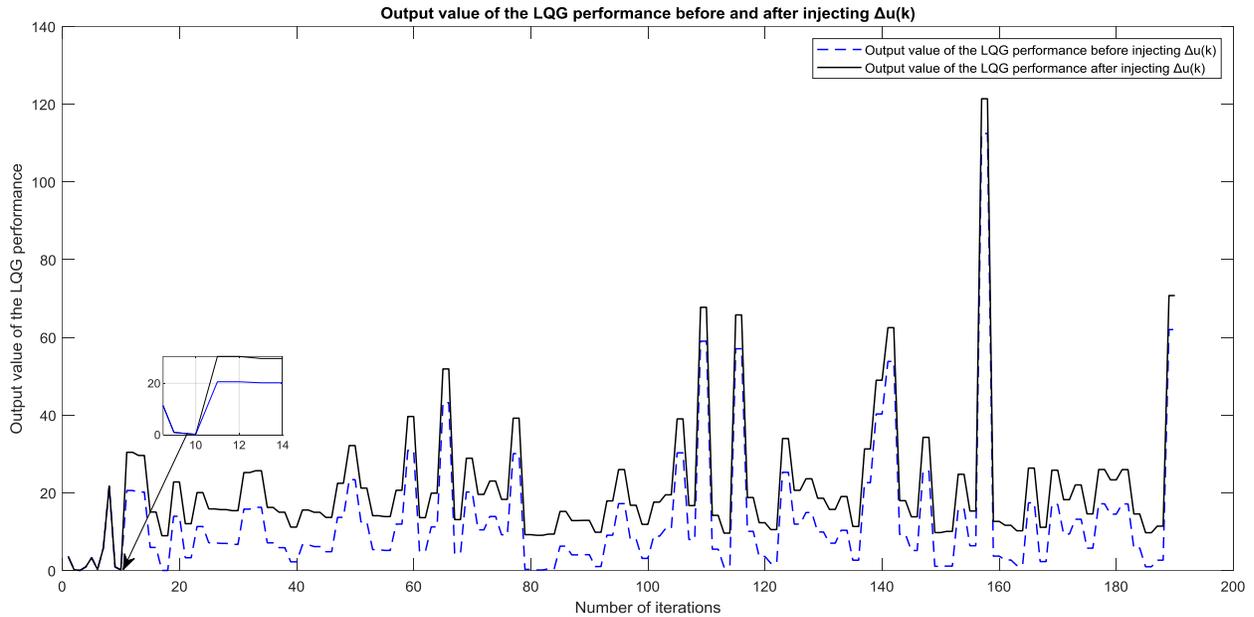

Figure 5

Although in previous approaches, the trade-off between attack detection delay or LQG performance was significant, in this approach based on the formulations and Figure 5, we can conclude that the difference in the mentioned trade-off is not considered when replay attack happens since the attack detection rate is rapid in first moments with the K-L divergence method rather than using $\chi^2$ detector or other statistical techniques, and therefore, defense strategies can be applied faster to the control system.

6. Conclusion

In this paper, a model of a replay attack is defined on CPSs and the performance of the control system under this kind of attack by proposing a new and different attack detection method, which is using the Kullback-Leibler divergence method is analyzed. It was proved that not for all control systems, the conventional estimation and attack detection methods are useful. Furthermore, to overcome this issue, for more general control systems with higher-order dynamics, a technique was provided to improve the detection rate considerably. In this case, the replay attack, which is the most frequent attack happening to control systems, can be stopped with defense strategies in the early moments with the proposed attack detection criterion.